\documentclass[12pt]{article}

\makeatletter
\renewcommand{\@oddfoot}{\hfill \thepage}
\makeatother

\usepackage[T2A]{fontenc}
\usepackage[cp1251]{inputenc}
\usepackage[russian, english]{babel}
\usepackage{mathrsfs}
\usepackage{amssymb}
\usepackage{amsmath,amsthm,amsfonts}
\usepackage{cite}
\usepackage[dvipdf]{graphicx}
\usepackage{fancyhdr}
\usepackage[titletoc]{appendix}

\begin{document}

\begin{center}
\Large\bf{Series representations for the characteristic function  
of the multidimensional Markov random flight}
\end{center}

\begin{center}
Alexander D. KOLESNIK\\
Institute of Mathematics and Computer Science\\
Academy Street 5, Kishinev 2028, Moldova\\
E-Mail: kolesnik@math.md
\end{center}

\begin{abstract} 
Two series representations of the characteristic function of the multidimensional symmetric Markov random flight with respect to Bessel functions and with respect to time variable, are given. Asymptotic formula for the second mixed moment function of the three-dimensional random flight is presented. 
\end{abstract}

\section{Introduction} 

\setcounter{equation}{0} 

Finite-velocity random motions in the Euclidean spaces of different dimensions whose evolution is driven by a homogeneous Poisson process, also called the Markov random flights, have become a very interesting and attractive model in recent decades. Historically, this field of stochastic analysis has started from the pioneering works by S. Goldstein \cite{gold} and M. Kac \cite{kac}, where a random motion of a particle moving with finite speed on the real line $\Bbb R^1$ and alternating its direction at Poissonian random time instants, was examined. The main result obtained in these works states that the transition density of the motion is the fundamental solution (the Green's function) to the hyperbolic telegraph equation (also called the damped wave equation), and can be found by solving it with respective initial conditions. Such one-dimensional stochastic motion, which is now referred to as the Goldstein-Kac telegraph process, as well as its numerous generalizations, have become the subject of a lot of researches, both theoretical and applied. The most important feature of the Goldstein-Kac telegraph process is that it can serve as an appropriate model of the one-dimensional transport that generates a diffusion with finite speed of propagation. For this reason, a generalized telegraph process has become a basis for describing many important physical properties of hyperbolic finite-velocity transport 
\cite{bras,giona1,giona2,giona3,giona4}. Relativistic aspects of the diffusion generated by finite-velocity random motions on the line, as well as their connections with some other physical processes, were given in \cite{cane1,cane2}. An approach to modeling the finite-velocity radioactive transfer based on some aspects of random flights theory was developed in \cite{eon}.  
A model of cosmic microwave background (CMB) radiation described by a telegraph equation with random initial conditions on the surface of the unit sphere in the three-dimensional space, as well as a number of its astrophysical interpretations, was presented in recent works \cite{broad1,broad2}. At present, the number of works devoted to telegraph processes and their applications in various fields of science, technology and engineering is really huge. The reader interested in the modern theory of telegraph processes and their applications in financial modeling may address to recent monograph \cite{kolrat}. 

The multidimensional counterpart of the Goldstein-Kac telegraph process is one of the most important 
generalizations performed by a particle moving at constant finite speed in the multidimensional Euclidean 
space $\Bbb R^m, \; m\ge 2$, and changing its direction at random Poissonian time instants by choosing the initial and new random directions on the surface of the unit $(m-1)$-dimensional sphere according to some probability law. Such multidimensional stochastic motions are now referred to as the Markov random flights. 
The main difference between the one-dimensional Goldstein-Kac telegraph process and its multidimensional counterparts steams from the fact that, while there are only two possible directions of motion on the real line $\Bbb R^1$, the multidimensional motion in the space $\Bbb R^m, \; m\ge 2$, has a continuum number of directions. Moreover, it turned out that the difficulty of analysis of the multidimensional Markov random flights fundamentally depends on the dimension of the phase space. While in the even-dimensional spaces 
$\Bbb R^2, \; \Bbb R^4$  and $\Bbb R^6$ one managed to obtain the transition densities of the motions in explicit forms (see \cite{mas,sta2,kol5,kol4,kol2}), in other dimensions, including the most interesting and important case of Markov random flight in the three-dimensional space $\Bbb R^3$, such analysis is turned out to be extremely difficult. For more details on the modern theory of multidimensional  Markov random flights, see recent monograph \cite{kol1}.

To overcome these difficulties, a method of integral transforms for describing the behaviour of multidimensional Markov random flights, was developed (see \cite{kol3} or \cite[Chapter 4]{kol1}). In the framework of this method, the basic characteristics of the Markov random flights in arbitrary dimensions 
in terms of Fourier-Laplace transforms, were given. One of the most important results is an explicit formula for the Fourier-Laplace transform of the transition density of the $m$-dimensional symmetric Markov random flight expressed in terms of Gauss hypergeometric function (see formula (\ref{eq01}) below). In order to get the most desirable goal, namely, the transition probability density, one needs to evaluate the double inverse Fourier and Laplace transforms of general relation (\ref{eq01}). However, this is an impracticable, at least directly, problem due to its very complicated form. That is why we make the first step in solving this problem and obtain two series representations of the characteristic function of the $m$-dimensional symmetric Markov random flight by evaluating the inverse Laplace transform of general relation (\ref{eq01}).

The article is organized as follows. In Section 2 we describe the $m$-dimensional symmetric Markov random flight and the structure of its distribution. We recall the general formula (\ref{eq01}) for the Fourier-Laplace transform of the transition density of the process. We also give a simple, but useful, auxiliary lemma concerning series decompositions of two special functions related to the Gauss hypergeometric function. This lemma will be applied in the next Sections 3 and 4 where two series representations of the characteristic function of the process with respect to Bessel functions and with respect to time variable 
$t$, will be obtained. The coefficients of these series representations are given by recurrent relations, 
as well as in the form of special determinants. Using these results in the final Section 5, we study the mixed moment function of the symmetric Markov random flight in the three-dimensional Euclidean space 
$\Bbb R^3$. We show that, for arbitrary time $t>0$, the first mixed moment is zero. For the second mixed moment function, we give an asymptotic relation with respect to the powers of time variable $t$. The error in this asymptotic relation has the order $o(t^7)$, that provides a very good approximation of the second mixed moment function for not too big values of time variable $t$.

\section{Preliminaries} 

\setcounter{equation}{0}

In this auxiliary section we give some basic notions and known results on the multidimensional 
symmetric Markov random flights which our further analysis will be relying on. 

Consider a particle that, at the initial time instant $t=0$, starts from the 
origin $\bold 0 = (0, \dots, 0)$ of the $m$-dimensional Euclidean space $\Bbb R^m, \; m\ge 2$. 
The particle moves with a constant finite speed $c$ (note that $c$ is treated as the 
constant norm of the velocity). The initial direction is a random $m$-dimensional vector 
uniformly distributed (the Lebesgue probability measure) on the 
unit sphere 
\begin{equation}\label{3.1.1}
S_1^m = \left\{ \bold x=(x_1, \dots ,x_m)\in \Bbb R^m: \;
\Vert\bold x\Vert^2 = x_1^2+ \dots +x_m^2=1 \right\} .
\end{equation} 
We emphasize that here and thereafter the upper index $m$ means the dimension of the space, 
which the sphere $S_1^m$ is considered in, not its own dimension which, clearly, is $m-1$. 
The particle changes its direction at random time instants that form a homogeneous Poisson 
flow of rate $\lambda>0$. In each of such Poissonian moments the particle instantly takes on a new random direction uniformly distributed on $S_1^m$ independently of its previous direction. 
Each sample path of this motion represents a broken line of total length $ct$ composed of 
the segments of $\lambda$-exponentially distributed random lengths and uniformly oriented 
in $\Bbb R^m$. The trajectories of the process are continuous and differentiable almost everywhere. 

Let $\bold X(t)=(X_1(t), \dots ,X_m(t))$ be the particle's position at 
arbitrary time instant $t>0$. The stochastic process $\bold X(t)$ is referred to as the 
{\it Markov random flight}. Let $N(t)$ denote the number of Poisson events that 
have occurred in the time interval $(0, t)$ and let $d\bold x$ be an infinitesimal 
element in the space $\Bbb R^m$ with Lebesgue measure $\mu(d\bold x) = dx_1 \dots dx_m$.

At arbitrary time moment $t>0$ the particle, with probability 1, is located in the 
$m$-dimensional ball of radius $ct$: 
\begin{equation}\label{3.1.2}
\bold B_{ct}^m = \left\{ \bold x=(x_1, \dots ,x_m)\in \Bbb R^m : \;
\Vert\bold x\Vert^2 = x_1^2+ \dots +x_m^2\le c^2t^2 \right\} .
\end{equation}

The distribution 
\begin{equation}\label{3.1.3}
\text{Pr} \left\{ \bold X(t)\in d\bold x \right\} = \text{Pr}
\left\{ X_1(t)\in dx_1, \dots , X_m(t)\in dx_m \right\} , \quad
\bold x\in\bold B_{ct}^m, \quad t>0, 
\end{equation}
consists of two components. The singular component of distribution (\ref{3.1.3}) 
is related to the case when no Poisson events occur in the time interval $(0,t)$ 
and, therefore, the particle does not change its initial direction. The singular component 
is concentrated on the sphere of radius $ct$: 
\begin{equation}\label{3.1.4}
S_{ct}^m =\partial\bold B_{ct}^m = \left\{ \bold x=(x_1,
\dots ,x_m)\in \Bbb R^m: \; \Vert\bold x\Vert^2 = x_1^2+ \dots
+x_m^2=c^2t^2 \right\} . 
\end{equation}
The probability of being on $S_{ct}^m$ at time moment $t$ is the probability 
that no Poisson events occur until time $t$ and it is equal to 
$$\text{Pr} \left\{ \bold X(t)\in S_{ct}^m \right\} = e^{-\lambda t} .$$
From this fact it follows that the density (in the sense of generalized functions) 
of the singular part of distribution (\ref{3.1.3}) has the form: 
\begin{equation}\label{3.1.5}
p^{(s)}(\bold x,t) = \frac{e^{-\lambda t}}{\text{mes} (S_{ct}^m)} \; \delta(c^2t^2 -
\Vert\bold x\Vert^2)  = \frac{e^{-\lambda t} \; \Gamma\left(
\frac{m}{2} \right)}{2\pi^{m/2} (ct)^{m-1}} \; \delta(c^2t^2 -
\Vert\bold x\Vert^2) , \qquad m\ge 2, 
\end{equation}
where $\delta(x)$ is the Dirac delta-function and $\text{mes}(S_{ct}^m)$ 
is the surface measure of sphere $S_{ct}^m$. 

If at least one Poisson event occurs in the time interval $(0,t)$ and, therefore, 
the particle at least once changes its direction, then, at time $t$, it is located 
strictly inside the ball $\bold B_{ct}^m$ and the probability of this event 
is equal to 
\begin{equation}\label{3.1.6}
\text{Pr} \left\{ \bold X(t)\in \text{int} \; \bold B_{ct}^m \right\} = 1-e^{-\lambda t} . 
\end{equation}
The part of distribution (\ref{3.1.3}) corresponding to this case is concentrated 
in the interior of the ball $\bold B_{ct}^m$ 
\begin{equation}\label{3.1.7}
\text{int} \; \bold B_{ct}^m = \left\{ \bold x=(x_1,
\dots ,x_m)\in \Bbb R^m: \; \Vert\bold x\Vert^2 = x_1^2+ \dots
+x_m^2<c^2t^2 \right\} 
\end{equation}
and forms its absolutely continuous component. Therefore, there exists the density of the absolutely continuous component of the distribution (\ref{3.1.3})
\begin{equation}\label{3.1.8}
p^{(ac)}(\bold x,t) = f(\bold x,t) \Theta(ct-\Vert\bold x\Vert) , \qquad 
\bold x\in \text{int} \; \bold B_{ct}^m , \quad t>0,
\end{equation}
where $f(\bold x,t)$ is some positive function absolutely continuous in $\text{int} \; \bold B_{ct}^m$ 
and $\Theta(x)$ is the Heaviside unit-step function. The existence of 
density (\ref{3.1.8}) follows from the fact that, since the sample paths of the process $\bold X(t)$ 
are continuous and differentiable almost everywhere, the distribution (\ref{3.1.3}) must contain an absolutely continuous component and this justifies the existence of density (\ref{3.1.8}). 
The existence of this density follows also from the fact that it can be represented as a Poissonian 
sum of convolutions. 

Hence, the density of distribution (\ref{3.1.3}) has the structure
$$p(\bold x,t) = p^{(s)}(\bold x,t) + p^{(ac)}(\bold x,t), 
\qquad \bold x\in\bold B_{ct}^m , \quad t>0,$$
where $p^{(s)}(\bold x,t)$ and $p^{(ac)}(\bold x,t)$ are the densities (in the sense 
of generalized functions) of the singular and absolutely continuous components 
of distribution (\ref{3.1.3}) given by (\ref{3.1.5}) and  (\ref{3.1.8}), respectively. 

Consider the characteristic function of the Markov random flight $\bold X(t)$ defined by  
\begin{equation}\label{3.2.3}
H(\boldsymbol\alpha,t) = \Bbb E \left\{ e^{i\langle\boldsymbol\alpha,\bold X(t)\rangle} \right\} ,  
\end{equation}
where $\boldsymbol\alpha=(\alpha_1, \dots ,\alpha_m) \in\Bbb R^m$ is the real-valued  
$m$-dimensional vector of inversion parameters and $\langle\boldsymbol\alpha,\bold X(t)\rangle$ 
is the inner product of the vectors $\boldsymbol\alpha$ and $\bold X(t)$. 

It is known (see \cite[Theorem 4.6.1]{kol1} or \cite[Formula (4.8)]{kol3}) that the Laplace transform 
of the characteristic function (\ref{3.2.3}) of the $m$-dimensional symmetric Markov random flight 
$\bold X(t)$  is given by the formula:  
\begin{equation}\label{eq01}
\mathcal L \left[ H(\boldsymbol\alpha,t) \right] (s) = \frac{F \left( \frac{1}{2},
\frac{m-2}{2}; \frac{m}{2};
\frac{(c\Vert\boldsymbol\alpha\Vert)^2}{(s+\lambda)^2
+(c\Vert\boldsymbol\alpha\Vert)^2} \right)}{
\sqrt{(s+\lambda)^2+(c\Vert\boldsymbol\alpha\Vert)^2} - \lambda \;
F \left( \frac{1}{2}, \frac{m-2}{2}; \frac{m}{2};
\frac{(c\Vert\boldsymbol\alpha\Vert)^2}{(s+\lambda)^2
+(c\Vert\boldsymbol\alpha\Vert)^2} \right)} , \quad m\ge 2,
\end{equation}
where  
\begin{equation}\label{eq02}
F(\xi,\eta;\zeta;z) \equiv \;
 _2F_1(\xi,\eta;\zeta;z)= \sum_{n=0}^{\infty} \dfrac{(\xi)_n
 (\eta)_n}{(\zeta)_n} \dfrac{z^n}{n!} 
\end{equation}
is the Gauss hypergeometric function and
$$(a)_n=a(a+1)\dots (a+n-1)=\frac{\Gamma (a+n)}{\Gamma (a)}$$
is the Pochhammer symbol. 

Function (\ref{eq01}) is holomorphic (analytical) 
and single-valued in the right half-plane of the complex plane, that is, for $\text{Re} \; s>0$.
It contains, in a hidden form, all the information concerning the distribution 
of $\bold X(t)$. In order to obtain the distribution, one needs to evaluate the inverse Laplace (with respect to complex parameter $s$) and Fourier (with respect to $m$-dimensional parameter 
$\boldsymbol\alpha$) transforms of function (\ref{eq01}). However, since it has a very complicated form, such double inversion seems to be an impracticable problem.  

In this article we make the first step in this direction and derive two series representations of the characteristic function of the process $\bold X(t)$ by evaluating the inverse Laplace transform of 
function (\ref{eq01}). The first representation is a functional series with respect to Bessel functions 
with variable indices. The second one is a series with respect to time variable $t$. 

Thus, in view of (\ref{eq01}) and the well-known properties of Laplace transform, 
we can write down the general formula for the characteristic function $H(t)$ as follows:
\begin{equation}\label{eq001}
H(\boldsymbol\alpha,t) =  e^{-\lambda t} \; \mathcal L_s^{-1} \left[ \frac{F \left( \frac{1}{2},
\frac{m-2}{2}; \frac{m}{2};
\frac{(c\Vert\boldsymbol\alpha\Vert)^2}{s^2
+(c\Vert\boldsymbol\alpha\Vert)^2} \right)}{
\sqrt{s^2+(c\Vert\boldsymbol\alpha\Vert)^2} - \lambda \;
F \left( \frac{1}{2}, \frac{m-2}{2}; \frac{m}{2};
\frac{(c\Vert\boldsymbol\alpha\Vert)^2}{s^2
+(c\Vert\boldsymbol\alpha\Vert)^2} \right)} \right] (t), \qquad m\ge 2,
\end{equation}
where $\mathcal L_s^{-1}$ means the inverse Laplace transformation with respect to inversion parameter $s$.
Our main efforts will be focused on evaluating the inverse Laplace transform in (\ref{eq001}).

The forthcoming analysis is substantially based on the following simple auxiliary lemma. 

{\bf Lemma 1.}\label{lemma1} 
{\it For arbitrary real constant $b$ and arbitrary positive constant $\lambda>0$, the following 
series decompositions hold:
\begin{equation}\label{eq1}
z F(\alpha, \beta; \gamma; b z^2) = \sum_{n=0}^{\infty} \xi_n \; z^n , 
\end{equation} 
\begin{equation}\label{eq2}
1-\lambda z F(\alpha, \beta; \gamma; b z^2) = \sum_{n=0}^{\infty} \eta_n \; z^n , 
\end{equation}
where the coefficients $\xi_n, \; \eta_n$ are given by the formulas: 
\begin{equation}\label{eq3}
\xi_n = \left\{ 
\aligned 0, \qquad & \quad \text{if} \;\; n=2r, \\
         \frac{(\alpha)_r \; (\beta)_r}{(\gamma)_r \;\; r!} \; b^r , & \quad \text{if} \;\; n=2r+1 , 
				\endaligned \right. \qquad r=0,1,2,\dots 
\end{equation}} 
\begin{equation}\label{eq4}
\eta_n = \left\{ 
\aligned 1, \qquad\qquad & \quad \text{if} \;\; n=0, \\
         0, \qquad\qquad & \quad \text{if} \;\; n=2r, \\
          - \lambda \; \frac{(\alpha)_{r-1} \; (\beta)_{r-1}}{(\gamma)_{r-1} \;\; (r-1)!} \; b^{r-1} , 
				  & \quad \text{if} \;\; n=2r-1 , 
				  \endaligned \right. \qquad r=1,2,3,\dots 
\end{equation} 

\begin{proof}
Using series representation of the Gauss hypergeometric function, we have: 
$$\aligned 
z F(\alpha, \beta; \gamma; b z^2) & = \sum_{l=0}^{\infty} 
\frac{(\alpha)_l \; (\beta)_l \; b^l}{(\gamma)_l \;\; l!} \; z^{2l+1} \\
& = \sum_{\substack{k=1\\k \; \text{is odd}}}^{\infty} 
\frac{(\alpha)_{(k-1)/2} \; (\beta)_{(k-1)/2} \;\; b^{(k-1)/2}}{(\gamma)_{(k-1)/2} \;\; 
\left(\frac{k-1}{2}\right)!} \; z^k \\
& = \sum_{n=0}^{\infty} \xi_n \; z^n , 
\endaligned$$
where the coefficients $\xi_n$ are given by (\ref{eq3}). 
Series representation (\ref{eq2}) immediately follows from (\ref{eq1}).
The lemma is proved. 
\end{proof}

In particular, formula (\ref{eq3}) yields: 
$$\xi_0=\xi_2=\xi_4=0, \qquad \xi_1=1, \quad \xi_3=b \; \frac{\alpha \beta}{\gamma}, \quad 
\xi_5 = \frac{b^2}{2} \; \frac{\alpha(\alpha+1) \beta(\beta+1)}{\gamma(\gamma+1)} .$$
Similarly, from (\ref{eq4}) we get: 
$$\eta_0=1, \quad \eta_1=-\lambda, \quad \eta_2=\eta_4=0, \quad  
\eta_3 = - \lambda b \; \frac{\alpha \beta}{\gamma} , \quad 
\eta_5 = - \frac{\lambda b^2}{2} \; \frac{\alpha(\alpha+1) \beta(\beta+1)}{\gamma(\gamma+1)} .$$

We see that the coefficients $\xi_n$ and $\eta_n$ are different for the index $n=0$, that is, $\xi_0=0$, 
while $\eta_0=1$. For all other indices, these coefficients are connected with each other by the relation 
\begin{equation}\label{eq013}
\eta_n = -\lambda \xi_n, \qquad n\ge 1.
\end{equation} 
Note also that both these coefficients are equal to zero for even indices, that is, 
$\xi_{2k}=\eta_{2k}=0, \; k\ge 1$.

\section{Series representation with respect to Bessel functions} 

In this section we derive a series representation of the characteristic function $H(\boldsymbol\alpha,t)$ 
of the $m$-dimensional symmetric Markov random flight $\bold X(t)$ with respect to Besel functions. Since the planar case ($m=2$) is well studied and its characteristic function is already known in the integral and series forms (see, for instance, \cite[Theorem 5.3.1]{kol1}), we omit this planar case and consider only 
the higher dimensions $m\ge 3$. The main result of this section is given by the following theorem. 

{\bf Theorem 1.}
{\it For arbitrary dimension $m\ge 3$, the characteristic function $H(\boldsymbol\alpha,t)$ of the 
$m$-dimensional symmetric Markov random flight $\bold X(t)$ has the following series representation: 
\begin{equation}\label{eq5}
H(\boldsymbol\alpha,t) = e^{-\lambda t} \; \sqrt{\pi} \; \sum_{n=1}^{\infty} 
\frac{\zeta_n(\boldsymbol\alpha)}{\Gamma\left(\frac{n}{2}\right)} 
\; \left(\frac{t}{2c\Vert\boldsymbol\alpha\Vert}\right)^{(n-1)/2} 
J_{(n-1)/2}(ct\Vert\boldsymbol\alpha\Vert) , 
\end{equation}
where $J_{\nu}(z)$ are Bessel functions and the coefficients $\zeta_n=\zeta_n(\boldsymbol\alpha)$ 
are given by the recurrent relation:
\begin{equation}\label{eq015} 
\zeta_1 = 1,  \qquad \zeta_n = \xi_n + \lambda \sum_{k=1}^{n-1} \zeta_{n-k} \; \xi_k , 
\qquad n\ge 2, 
\end{equation}
with 
\begin{equation}\label{eq012}
\xi_n = \left\{ 
\aligned 0, \hskip 3cm & \quad \text{if} \;\; n=2r, \\
         \frac{(2r-1)!!}{(2r)!!} \; \frac{m-2}{2r+m-2} \; 
				 (c\Vert\boldsymbol\alpha\Vert)^{2r} , 
				 & \quad \text{if} \;\; n=2r+1 , 
				 \endaligned \right. \qquad r=0,1,2,\dots 
\end{equation}
where 
$(2r)!! = 2\cdot 4\cdot 6\cdot \; \dots \; \cdot (2r), \;\; 
(2r-1)!! = 1\cdot 3\cdot 5\cdot \; \dots\; \cdot (2r-1), \; 
r\ge 0, \;\; 0!!=(-1)!! \overset{\text{def}}{=} 1$.} 

\begin{proof}
Consider separately the fraction in square brackets in (\ref{eq001}). 
Dividing the numerator and denominator of this fraction by 
$\sqrt{s^2+(c\Vert\boldsymbol\alpha\Vert)^2}$ and introducing the variable 
\begin{equation}\label{eq6}
z = \frac{1}{\sqrt{s^2+(c\Vert\boldsymbol\alpha\Vert)^2}}
\end{equation}
we can represent it as follows: 
\begin{equation}\label{eq7}
A(z) = \frac{z \; F \left( \frac{1}{2}, \frac{m-2}{2}; \frac{m}{2};
c^2\Vert\boldsymbol\alpha\Vert^2 z^2\right)}{
1 - \lambda z \; F \left( \frac{1}{2}, \frac{m-2}{2}; \frac{m}{2};
c^2\Vert\boldsymbol\alpha\Vert^2 z^2\right)}
\end{equation} 
Applying Lemma 1 to this fraction, we get: 
\begin{equation}\label{eq8}
A(z) = \frac{\sum\limits_{n=0}^{\infty} \xi_n \; z^n }{\sum\limits_{n=0}^{\infty} \eta_n \; z^n} ,
\end{equation} 
where, according to (\ref{eq3}) and (\ref{eq4}), the coefficients $\xi_n, \; \eta_n$ are given by 
\begin{equation}\label{eq9}
\xi_n = \left\{ 
\aligned 0, \qquad & \quad \text{if} \;\; n=2r, \\
         \frac{\left(\frac{1}{2}\right)_r \; \left(\frac{m-2}{2}\right)_r \; 
				(c\Vert\boldsymbol\alpha\Vert)^{2r}}{\left(\frac{m}{2}\right)_r \;\; r!} , 
				& \quad \text{if} \;\; n=2r+1 , 
				\endaligned \right. \qquad r=0,1,2,\dots 
\end{equation} 
\vskip 0.2cm
\begin{equation}\label{eq10}
\eta_n = \left\{ 
\aligned 1, \qquad\qquad & \quad \text{if} \;\; n=0, \\
         0, \qquad\qquad & \quad \text{if} \;\; n=2r, \\
          - \lambda \; \frac{\left(\frac{1}{2}\right)_{r-1} \; \left(\frac{m-2}{2}\right)_{r-1} \; 
					(c\Vert\boldsymbol\alpha\Vert)^{2(r-1)}}{\left(\frac{m}{2}\right)_{r-1} \;\; (r-1)!} , 
				  & \quad \text{if} \;\; n=2r-1 , 
				  \endaligned \right. \qquad r=1,2,3,\dots 
\end{equation} 
Coefficients (\ref{eq9}) and (\ref{eq10}) depend on the vector variable $\boldsymbol\alpha$, but, for the sake of brevity, we omit it thereafter, as well as in the coefficients 
$\zeta_n=\zeta_n(\boldsymbol\alpha)$ (see recurrent relation (\ref{eq15}) below), bearing in mind, however, that all these coefficients are, in fact, functions of the vector variable $\boldsymbol\alpha$ (more precisely, of its norm $\Vert\boldsymbol\alpha\Vert$). 

Using the definition of Pochhammer symbol and the relation 
$$\frac{(a)_k}{(a+1)_k} = \frac{a}{a+k} ,$$
as well as the well-known relations for Euler gamma-function: 
\begin{equation}\label{eq11} 
\Gamma\left(\frac{1}{2}\right) = \sqrt{\pi} , \qquad \Gamma\left(k+\frac{1}{2}\right) = 
\frac{\sqrt{\pi}}{2^k} (2k-1)!! , \qquad \Gamma(x+1) = x \Gamma(x) , \qquad  k\in\Bbb N, 
\end{equation} 
one can easily show that 
$$\left(\frac{1}{2}\right)_r = 2^{-2r} \; \frac{(2r)!}{r!} , \qquad 
\left(\frac{1}{2}\right)_{r-1} = 2^{-2(r-1)} \; \frac{(2r-2)!}{(r-1)!} ,$$ 
$$\frac{\left(\frac{m-2}{2}\right)_r}{\left(\frac{m}{2}\right)_r} = \frac{m-2}{2r+m-2} , \qquad 
\frac{\left(\frac{m-2}{2}\right)_{r-1}}{\left(\frac{m}{2}\right)_{r-1}} = \frac{m-2}{2r+m-4} ,$$
and, therefore, coefficients (\ref{eq9}) and (\ref{eq10}) take the form: 
\begin{equation}\label{eq12}
\xi_n = \left\{ 
\aligned 0, \hskip 3cm & \quad \text{if} \;\; n=2r, \\
         \frac{(2r-1)!!}{(2r)!!} \; \frac{m-2}{2r+m-2} \; 
				 (c\Vert\boldsymbol\alpha\Vert)^{2r} , 
				 & \quad \text{if} \;\; n=2r+1 , 
				 \endaligned \right. \qquad r=0,1,2,\dots 
\end{equation} 

\vskip 0.2cm

\begin{equation}\label{eq13}
\eta_n = \left\{ 
\aligned 1, \hskip 3cm & \quad \text{if} \;\; n=0, \\
         0, \hskip 3cm & \quad \text{if} \;\; n=2r, \\
          - \lambda \; \frac{(2r-3)!!}{(2r-2)!!} \; \frac{m-2}{2r+m-4} \; 
					(c\Vert\boldsymbol\alpha\Vert)^{2r-2},  
				  & \quad \text{if} \;\; n=2r-1 , 
				  \endaligned \right. \qquad r=1,2,3,\dots 
\end{equation} 
where 
$(2r)!! = 2\cdot 4\cdot 6\cdot \; \dots \; \cdot (2r), \;\; 
(2r-1)!! = 1\cdot 3\cdot 5\cdot \; \dots\; \cdot (2r-1), \; 
r\ge 0, \;\; 0!!=(-1)!! \overset{\text{def}}{=} 1$. 

Applying now a formula for the quotient of two power series (see, for instance, \cite[Formula 0.313]{gr} 
or \cite[page 754, item 4]{pbm1}) to fraction (\ref{eq8}) and taking into account (\ref{eq6}), we get: 
\begin{equation}\label{eq14}
A(s) = \sum_{n=0}^{\infty} \frac{\zeta_n}{(s^2+(c\Vert\boldsymbol\alpha\Vert)^2)^{n/2}} ,
\end{equation} 
where the coefficients $\zeta_n, \; n\ge 0,$ are given by the recurrent relation:  
\begin{equation}\label{eq15} 
\zeta_0 = 0, \qquad \zeta_1 = 1,  \qquad \zeta_n = \xi_n - \sum_{k=1}^{n-1} \zeta_{n-k} \; \eta_k , 
\qquad n\ge 2. 
\end{equation} 
In particular, in view of (\ref{eq12}) and (\ref{eq13}), we have: 
\begin{equation}\label{eq16} 
\aligned  
\zeta_0 & = 0, \qquad \zeta_1 = 1, \quad \zeta_2 = \lambda, \qquad \zeta_3 = \lambda^2 + 
\frac{m-2}{2m} \; (c\Vert\boldsymbol\alpha\Vert)^2 , \\ \\
\zeta_4 & = \lambda^3 + \lambda \; \frac{m-2}{m} \; (c\Vert\boldsymbol\alpha\Vert)^2   
\qquad \zeta_5 = \lambda^4 + \frac{3}{2} \; \lambda^2 \; \frac{m-2}{m} \; (c\Vert\boldsymbol\alpha\Vert)^2 
+ \frac{3}{8} \; \frac{m-2}{m+2} \; (c\Vert\boldsymbol\alpha\Vert)^4 .
\endaligned
\end{equation} 

Series (\ref{eq14}) is uniformly convergent in the right half-plane of the complex plane. Applying 
the inverse Laplace transformation $\mathcal L_s^{-1}$ to this series (\ref{eq14}) and taking into account 
(see, for instance \cite[Table 8.4-1, Formula 57]{korn}) that 
$$\mathcal L_s^{-1} \left[ \frac{1}{(s^2+a^2)^k} \right](t) = 
\frac{\sqrt{\pi}}{\Gamma(k)} \; \left( \frac{t}{2a} \right)^{k-1/2} J_{k-1/2}(at) , \qquad k>0,$$
where $J_{\nu}(z)$ is Bessel function, we obtain: 
$$\aligned 
\mathcal L_s^{-1} [ A(s)](t) & = \sum_{n=1}^{\infty} \zeta_n \; 
\mathcal L_s^{-1} \left[ \frac{1}{(s^2+(c\Vert\boldsymbol\alpha\Vert)^2)^{n/2}} \right](t) \\ 
& = \sqrt{\pi} \; \sum_{n=1}^{\infty} \frac{\zeta_n}{\Gamma\left(\frac{n}{2}\right)} \; 
\left(\frac{t}{2c\Vert\boldsymbol\alpha\Vert}\right)^{(n-1)/2} J_{(n-1)/2}(ct\Vert\boldsymbol\alpha\Vert) .
\endaligned$$
Returning to (\ref{eq001}), we arrive at the statement of the theorem. 
\end{proof}

\bigskip 

{\bf Remark 1.} In view of (\ref{eq16}) and taking into account that  
\begin{equation}\label{eq17} 
\aligned 
& \Gamma\left( \frac{1}{2} \right) = \sqrt{\pi}, \quad \Gamma\left( \frac{3}{2} \right) = 
\frac{\sqrt{\pi}}{2}, \quad \Gamma\left( \frac{5}{2} \right) = \frac{3\sqrt{\pi}}{4}, \\ 
& J_{1/2}(z) = \sqrt{\frac{2}{\pi z}} \; \sin{z} , \quad J_{3/2}(z) = \sqrt{\frac{2}{\pi z}} 
\left( \frac{\sin{z}}{z} - \cos{z} \right) ,
\endaligned
\end{equation} 
we can write down the five terms of series (\ref{eq5}):
\begin{equation}\label{eq18} 
\aligned 
H(\boldsymbol\alpha,t) & = e^{-\lambda t} \biggl[ J_0(ct\Vert\boldsymbol\alpha\Vert) + \lambda \; 
\frac{\sin(ct\Vert\boldsymbol\alpha\Vert)}{c\Vert\boldsymbol\alpha\Vert} + 
\left( \lambda^2 + \frac{m-2}{2m} \; (c\Vert\boldsymbol\alpha\Vert)^2 \right)  
\frac{t}{c\Vert\boldsymbol\alpha\Vert} \; J_1(ct\Vert\boldsymbol\alpha\Vert) \\ 
& + \left( \lambda^3 + \lambda \; \frac{m-2}{m} \; (c\Vert\boldsymbol\alpha\Vert)^2 \right)  
\frac{t}{2(c\Vert\boldsymbol\alpha\Vert)^2} \; 
\left( \frac{\sin(ct\Vert\boldsymbol\alpha\Vert)}{ct\Vert\boldsymbol\alpha\Vert} - 
\cos(ct\Vert\boldsymbol\alpha\Vert) \right) \\ 
& + \left( \lambda^4 + \frac{3}{2} \; \lambda^2 \; \frac{m-2}{m} \; 
(c\Vert\boldsymbol\alpha\Vert)^2 + \frac{3}{8} \; \frac{m-2}{m+2} \;  
(c\Vert\boldsymbol\alpha\Vert)^4 \right)  
\frac{t^2}{3(c\Vert\boldsymbol\alpha\Vert)^2} \; 
J_2(ct\Vert\boldsymbol\alpha\Vert) + \dots \biggr] .
\endaligned
\end{equation} 

\bigskip 

{\bf Remark 2.} In view of \cite[Formula 0.313]{gr} or \cite[page 754, item 4]{pbm1}, relation 
(\ref{eq013}) and taking into account that $\xi_1=1, \; \xi_n=0, \; n=2r, \; r=0,1,2,\dots,$ and 
$\eta_0=1, \; \eta_1=-\lambda, \; \eta_n=0, \; n=2r, \; r=1,2,\dots,$ we can, instead of recurrent 
relation (\ref{eq015}), give an explicit formula for evaluating the coefficients 
$\zeta_n, \; n\ge 2,$ in the form of the following $(n\times n)$-determinant: 
\begin{equation}\label{eq19} 
\zeta_n = (-1)^{n+1}   
\begin{vmatrix}
1 & 1 & 0 & 0 & 0 & 0 \; \dots \; 0 & 0 \\
0 & -\lambda & 1 & 0 & 0 & 0 \; \dots \; 0 & 0 \\ 
\xi_3 & 0 & -\lambda & 1 & 0 & 0 \; \dots \; 0 & 0\\ 
0 & -\lambda\xi_3 & 0 & -\lambda & 1 & 0 \; \dots \; 0 & 0 \\ 
\xi_5 & 0 & -\lambda\xi_3 & 0 & -\lambda & 1 \; \dots \; 0 & 0\\  
\vdots &  \vdots & \vdots & \vdots & \vdots & \vdots \qquad \; \vdots & \vdots\\  
\xi_{n-1} & -\lambda\xi_{n-2} & -\lambda\xi_{n-3} & -\lambda\xi_{n-4} & -\lambda\xi_{n-5} & \;\;\; 0 \; \dots\;  -\lambda & 1 \\ 
\xi_n & -\lambda\xi_{n-1} & -\lambda\xi_{n-2} & -\lambda\xi_{n-3} & -\lambda\xi_{n-4} & -\lambda\xi_{n-5} \; \dots\; 0 & -\lambda  
\end{vmatrix} , 
\qquad n\ge 2. 
\end{equation} 
In particular, 
$$\zeta_2 = -  
\begin{vmatrix} 
1 & 1\\ 
0 & -\lambda
\end{vmatrix} 
= \lambda , \qquad 
\zeta_3 = 
\begin{vmatrix} 
1 & 1 & 0\\ 
0 & -\lambda & 1\\ 
\xi_3 & 0 & -\lambda  
\end{vmatrix} 
= \lambda^2 + \xi_3 = \lambda^2 + \frac{m-2}{2m} \; (c\Vert\boldsymbol\alpha\Vert)^2 ,$$ 
$$\aligned 
\zeta_4 = & - 
\begin{vmatrix} 
1 & 1 & 0 & 0 \\
0 & -\lambda & 1 & 0\\ 
\xi_3 & 0 & -\lambda & 1\\
0 & -\lambda\xi_3 & 0 & -\lambda 
\end{vmatrix} 
= - 
\begin{vmatrix}  
-\lambda & 1 & 0\\ 
0 & -\lambda & 1\\
-\lambda\xi_3 & 0 & -\lambda 
\end{vmatrix} + 
\begin{vmatrix} 
0 & 1 & 0\\ 
\xi_3 & -\lambda & 1\\
0 & 0 & -\lambda 
\end{vmatrix} \\ 
\\ 
= & \lambda^3 + 2\lambda \xi_3 = 
\lambda^3 + \lambda \; \frac{m-2}{m} \; (c\Vert\boldsymbol\alpha\Vert)^2  , 
\endaligned$$ 
\vskip 0.2cm
$$\aligned 
\zeta_5 & =  
\begin{vmatrix} 
1 & 1 & 0 & 0 & 0\\
0 & -\lambda & 1 & 0 & 0\\ 
\xi_3 & 0 & -\lambda & 1 & 0\\
0 & -\lambda\xi_3 & 0 & -\lambda & 1\\
\xi_5 & 0 & -\lambda\xi_3 & 0 & -\lambda   
\end{vmatrix} 
=  
\begin{vmatrix} 
-\lambda & 1 & 0 & 0\\ 
0 & -\lambda & 1 & 0\\
-\lambda\xi_3 & 0 & -\lambda & 1\\
0 & -\lambda\xi_3 & 0 & -\lambda   
\end{vmatrix} 
- 
\begin{vmatrix} 
0 & 1 & 0 & 0\\ 
\xi_3 & -\lambda & 1 & 0\\
0 & 0 & -\lambda & 1\\
\xi_5 & -\lambda\xi_3 & 0 & -\lambda   
\end{vmatrix} \\
& = -\lambda 
\begin{vmatrix} 
-\lambda & 1 & 0\\
0 & -\lambda & 1\\
-\lambda\xi_3 & 0 & -\lambda   
\end{vmatrix} 
- 
\begin{vmatrix} 
0 & 1 & 0\\
-\lambda\xi_3 & -\lambda & 1\\
0 & 0 & -\lambda   
\end{vmatrix} 
+ 
\begin{vmatrix} 
\xi_3 & 1 & 0\\
0 & -\lambda & 1\\
\xi_5 & 0 & -\lambda   
\end{vmatrix} \\ 
\\
& = \lambda^4 + 3\lambda^2 \xi_3 + \xi_5 =  
\lambda^4 + \frac{3}{2} \; \lambda^2 \; \frac{m-2}{m} \; (c\Vert\boldsymbol\alpha\Vert)^2 
+ \frac{3}{8} \; \frac{m-2}{m+2} \; (c\Vert\boldsymbol\alpha\Vert)^4 ,
\endaligned$$ 
and this exactly coincides with (\ref{eq16}).

\section{Series representation with respect to time variable} 

Theorem 1 yields a series representation of the characteristic function $H(\boldsymbol\alpha,t)$ with respect to Bessel functions with variable indices. However, one can also obtain a sereis representation 
of $H(\boldsymbol\alpha,t)$ with respect to powers of time variable $t$. This result is given by the following theorem. 

{\bf Theorem 2.} {\it For arbitrary dimension $m\ge 3$, the characteristic function $H(\boldsymbol\alpha,t)$of the $m$-dimensional symmetric Markov random flight $\bold X(t)$  has the following series representation: 
\begin{equation}\label{eq20} 
H(\boldsymbol\alpha,t) = e^{-\lambda t} \sum_{n=0}^{\infty} 
\frac{\gamma_{n+1}(\boldsymbol\alpha)}{n!} \; t^n ,
\end{equation} 
where the coefficients $\gamma_n=\gamma_n(\boldsymbol\alpha)$ are given by the recurrent relation: 
\begin{equation}\label{eq024} 
\gamma_1 = 1,  \qquad \gamma_n = \theta_n + \lambda \sum_{k=1}^{n-1} 
\gamma_{n-k} \; \theta_k , \qquad n\ge 2, 
\end{equation} 
and  
\begin{equation}\label{eq0090}
\theta_n = \left\{ 
\aligned 0, \qquad & \quad \text{if} \;\; n=2r, \\
         (-1)^r \; \frac{\left( \frac{1}{2} \right)_r}{\left( \frac{m}{2} \right)_r} \; 
				(c\Vert\boldsymbol\alpha\Vert)^{2r} , 
				& \quad \text{if} \;\; n=2r+1 , 
				\endaligned \right. \qquad r=0,1,2,\dots 
\end{equation}}

\begin{proof} 
According to \cite[Formula 9.131(1)]{gr},  
$$\frac{1}{\sqrt{s^2+(c\Vert\boldsymbol\alpha\Vert)^2}} \;
F \left( \frac{1}{2}, \frac{m-2}{2}; \frac{m}{2};
\frac{(c\Vert\boldsymbol\alpha\Vert)^2}{s^2+(c\Vert\boldsymbol\alpha\Vert)^2}
\right) = \frac{1}{s} \; F \left( \frac{1}{2}, 1; \frac{m}{2};
-\frac{(c\Vert\boldsymbol\alpha\Vert)^2}{s^2} \right) ,$$
and, therefore, by introducing the variable 
\begin{equation}\label{eq21}  
z = \frac{1}{s}
\end{equation}
we can rewrite formula (\ref{eq7}) as follows: 
\begin{equation}\label{eq22}  
A(z) = \frac{z \; F \left( \frac{1}{2}, 1; \frac{m}{2};
-c^2\Vert\boldsymbol\alpha\Vert^2 z^2\right)}{
1 - \lambda z \; F \left( \frac{1}{2}, 1; \frac{m}{2};
-c^2\Vert\boldsymbol\alpha\Vert^2 z^2\right)} .
\end{equation}
Applying Lemma 1 to this fraction, we get:  
\begin{equation}\label{eq08}
A(z) = \frac{\sum\limits_{n=0}^{\infty} \theta_n \; z^n }{\sum\limits_{n=0}^{\infty} \sigma_n \; z^n} ,
\end{equation} 
where the coefficients $\theta_n, \; \sigma_n$ are given by 
\begin{equation}\label{eq09}
\theta_n = \left\{ 
\aligned 0, \qquad & \quad \text{if} \;\; n=2r, \\
         (-1)^r \; \frac{\left(\frac{1}{2}\right)_r \; (1)_r \; 
				(c\Vert\boldsymbol\alpha\Vert)^{2r}}{\left(\frac{m}{2}\right)_r \;\; r!} , 
				& \quad \text{if} \;\; n=2r+1 , 
				\endaligned \right. \qquad r=0,1,2,\dots 
\end{equation} 
\vskip 0.2cm
\begin{equation}\label{eq010}
\sigma_n = \left\{ 
\aligned 1, \qquad\qquad & \quad \text{if} \;\; n=0, \\
         0, \qquad\qquad & \quad \text{if} \;\; n=2r, \\
          (-1)^r \; \lambda \; \frac{\left(\frac{1}{2}\right)_{r-1} \; (1)_{r-1} \; 
					(c\Vert\boldsymbol\alpha\Vert)^{2(r-1)}}{\left(\frac{m}{2}\right)_{r-1} \;\; (r-1)!} , 
				  & \quad \text{if} \;\; n=2r-1 , 
				  \endaligned \right. \qquad r=1,2,3,\dots 
\end{equation} 
Coefficients (\ref{eq09}) and (\ref{eq010}) depend on the vector variable $\boldsymbol\alpha$, but, for the sake of brevity, we omit it thereafter, as well as in the coefficients 
$\gamma_n=\gamma_n(\boldsymbol\alpha)$ (see recurrent relation (\ref{eq0015}) below), bearing in mind, however, that all these coefficients are, in fact, functions of the vector variable $\boldsymbol\alpha$.   

Since $(1)_r = r!, \; (1)_{r-1} = (r-1)!$, the coefficients (\ref{eq09}) and (\ref{eq010}) reduce to  
\begin{equation}\label{eq009}
\theta_n = \left\{ 
\aligned 0, \qquad & \quad \text{if} \;\; n=2r, \\
         (-1)^r \; \frac{\left( \frac{1}{2} \right)_r}{\left( \frac{m}{2} \right)_r} \; 	
				 (c\Vert\boldsymbol\alpha\Vert)^{2r} , 
				 & \quad \text{if} \;\; n=2r+1 , 
				 \endaligned \right. \qquad r=0,1,2,\dots 
\end{equation} 
\vskip 0.2cm
\begin{equation}\label{eq0100}
\sigma_n = \left\{ 
\aligned 1, \qquad\qquad & \quad \text{if} \;\; n=0, \\
         0, \qquad\qquad & \quad \text{if} \;\; n=2r, \\
          (-1)^r \;  \lambda \; \frac{\left( \frac{1}{2} \right)_{r-1}}{\left( \frac{m}{2} \right)_{r-1}} 
					\; (c\Vert\boldsymbol\alpha\Vert)^{2(r-1)} ,  
				  & \quad \text{if} \;\; n=2r-1 , 
				  \endaligned \right. \qquad r=1,2,3,\dots 
\end{equation}
In particular, 
\begin{equation}\label{eq23} 
\aligned 
& \theta_0 = \theta_2 = \theta_4 = \theta_6 = 0, \qquad 
\theta_1=1, \quad \theta_3 = - \frac{1}{m} \; (c\Vert\boldsymbol\alpha\Vert)^2, \\ 
& \theta_5 = \frac{3!!}{m(m+2)} \; (c\Vert\boldsymbol\alpha\Vert)^4, \quad 
\theta_7 = - \frac{5!!}{m (m+2) (m+4)} \; (c\Vert\boldsymbol\alpha\Vert)^6, \\
& \sigma_0=1, \quad \sigma_2=\sigma_4=0, \quad \sigma_1=-\lambda, \quad 
\sigma_3 = \lambda \; \frac{1}{m} \; (c\Vert\boldsymbol\alpha\Vert)^2, \quad 
\sigma_5 = -\lambda \; \frac{3}{m(m+2)} \; (c\Vert\boldsymbol\alpha\Vert)^4 . 
\endaligned
\end{equation}

Applying now a formula for the quotient of two power series (see, for instance, \cite[Formula 0.313]{gr} 
or \cite[page 754, item 4]{pbm1}) to fraction (\ref{eq08}) and taking into account (\ref{eq21}), we get: 
\begin{equation}\label{eq014}
A(s) = \sum_{n=0}^{\infty} \frac{\gamma_n}{s^n} ,
\end{equation} 
where the coefficients $\gamma_n, \; n\ge 0,$ are given by the recurrent relation: 
\begin{equation}\label{eq0015} 
\gamma_0 = 0, \qquad \gamma_1 = 1,  \qquad \gamma_n = \theta_n - \sum_{k=1}^{n-1} \gamma_{n-k} \; \sigma_k , 
\qquad n\ge 2. 
\end{equation} 
Taking onto account that $\sigma_n=-\lambda\theta_n, \; n\ge 1,$ we can rewrite recurrent relation 
(\ref{eq0015}) as follows: 
\begin{equation}\label{eq24} 
\gamma_0 = 0, \qquad \gamma_1 = 1,  \qquad \gamma_n = \theta_n + \lambda \sum_{k=1}^{n-1} 
\gamma_{n-k} \; \theta_k , \qquad n\ge 2. 
\end{equation}
In particular, some simple calculations yield: 
\begin{equation}\label{eq25} 
\aligned 
& \gamma_0 = 0, \qquad \gamma_1 = 1, \qquad \gamma_2 = \lambda, \qquad 
\gamma_3 = \lambda^2 - \frac{1}{m} \; (c\Vert\boldsymbol\alpha\Vert)^2 , \\  
& \gamma_4 = \lambda^3 - \lambda \; \frac{2}{m} \; (c\Vert\boldsymbol\alpha\Vert)^2 , \qquad 
\gamma_5 = \lambda^4 - \lambda^2 \; \frac{3}{m} \; (c\Vert\boldsymbol\alpha\Vert)^2 
+ \frac{3}{m(m+2)} \; (c\Vert\boldsymbol\alpha\Vert)^4 \\
& \gamma_6 = \lambda^5 - \lambda^3 \; \frac{4}{m} \; (c\Vert\boldsymbol\alpha\Vert)^2 
+ \lambda \; \frac{7m+2}{m^2 (m+2)} \; (c\Vert\boldsymbol\alpha\Vert)^4 \\ 
& \gamma_7 = \lambda^6 - \lambda^4 \; \frac{5}{m} \; (c\Vert\boldsymbol\alpha\Vert)^2 
+ \lambda^2 \; \frac{12m+6}{m^2 (m+2)} \; (c\Vert\boldsymbol\alpha\Vert)^4  
- \frac{15}{m(m+2)(m+4)} \; (c\Vert\boldsymbol\alpha\Vert)^6 \\ 
& \gamma_8 = \lambda^7 - \lambda^5 \; \frac{6}{m} \; (c\Vert\boldsymbol\alpha\Vert)^2  
+ \lambda^3 \; \frac{18m+12}{m^2 (m+2)} \; (c\Vert\boldsymbol\alpha\Vert)^4  
- \lambda \; \frac{36m+24}{m^2 (m+2)(m+4)} \; (c\Vert\boldsymbol\alpha\Vert)^6 \\ 
\endaligned 
\end{equation}

Function $A(s)$ given by series (\ref{eq014}) is holomorphic (analytical) everywhere in the right half-plane 
of the complex plane, that is, for $\text{Re} \; s>0$. Therefore, by applying the inverse Laplace transformation $\mathcal L_s^{-1}$ to series (\ref{eq014}) and taking 
into account (see, for instance \cite[Table 8.4-1, Formula 3]{korn}) that 
$$\mathcal L_s^{-1} \left[ \frac{1}{s^n} \right](t) = \frac{t^{n-1}}{(n-1)!} , \qquad n=1,2,\dots ,$$ 
and that $\gamma_0=0$, we obtain: 
$$\mathcal L_s^{-1} [A(s)](t) = \sum_{n=1}^{\infty} \gamma_n \; \frac{t^{n-1}}{(n-1)!} = 
\sum_{n=0}^{\infty} \frac{\gamma_{n+1}}{n!} \; t^n .$$
Multiplying this by $e^{-\lambda t}$, we arrive at the statement of the theorem. 
\end{proof}

\bigskip 

{\bf Remark 3.} Taking into account (\ref{eq25}), we can write down the five terms of 
the characteristic function (\ref{eq20}):   
\begin{equation}\label{eq26} 
\aligned 
H(\boldsymbol\alpha,t) = e^{-\lambda t} \biggl[ & 1 + \lambda t + \frac{1}{2!} 
\left( \lambda^2 - \frac{1}{m} \; (c\Vert\boldsymbol\alpha\Vert)^2  \right) t^2  
 + \frac{1}{3!} \left( \lambda^3 - \lambda \; \frac{2}{m} \; (c\Vert\boldsymbol\alpha\Vert)^2 \right) t^3 \\
& + \frac{1}{4!} \left( \lambda^4 - \lambda^2 \; \frac{3}{m} \; (c\Vert\boldsymbol\alpha\Vert)^2 
+ \frac{3}{m(m+2)} \; (c\Vert\boldsymbol\alpha\Vert)^4 \right) t^4 + \dots \biggr] .
\endaligned 
\end{equation}

\bigskip 

{\bf Remark 4.} Similarly as above (see Remark 2), coefficients $\gamma_n$ can be represented in the explicit form of the following $(n\times n)$-determinant: 
\begin{equation}\label{eq27} 
\gamma_n = (-1)^{n+1}   
\begin{vmatrix}
1 & 1 & 0 & 0 & 0 & 0 \; \dots \; 0 & 0 \\
0 & -\lambda & 1 & 0 & 0 & 0 \; \dots \; 0 & 0 \\ 
\theta_3 & 0 & -\lambda & 1 & 0 & 0 \; \dots \; 0 & 0\\ 
0 & -\lambda\theta_3 & 0 & -\lambda & 1 & 0 \; \dots \; 0 & 0 \\ 
\theta_5 & 0 & -\lambda\theta_3 & 0 & -\lambda & 1 \; \dots \; 0 & 0\\  
\vdots &  \vdots & \vdots & \vdots & \vdots & \vdots \qquad \; \vdots & \vdots\\  
\theta_{n-1} & -\lambda\theta_{n-2} & -\lambda\theta_{n-3} & -\lambda\theta_{n-4} & -\lambda\theta_{n-5} & \;\;\; 0 \; \dots\;  -\lambda & 1 \\ 
\theta_n & -\lambda\theta_{n-1} & -\lambda\theta_{n-2} & -\lambda\theta_{n-3} & -\lambda\theta_{n-4} & -\lambda\theta_{n-5} \; \dots\; 0 & -\lambda  
\end{vmatrix} , 
\qquad n\ge 2. 
\end{equation} 
Similarly to Remark 2, one can easily check that, for $n=2, 3, 4, 5$, determinant (\ref{eq27}) produces the same coefficients (\ref{eq25}). 

\bigskip 

{\bf Remark 5.} One can see that, for arbitrary $n\ge 1$, the coefficient 
$\gamma_n = \gamma_n(\Vert\boldsymbol\alpha\Vert)$ is a polynomial of the variable 
$\Vert\boldsymbol\alpha\Vert$ of the power $2\left[ \frac{n-1}{2} \right], \; n\ge 1,$ or, equivalently, 
a polynomial of the variable $(\alpha_1^2+ \dots +\alpha_m^2)$ of the power 
$\left[ \frac{n-1}{2} \right], \; n\ge 1$, where $[ \; \cdot \; ]$ means the integer part of a number. 

\bigskip 

{\bf Remark 6.} The forms of coefficients $\gamma_n$ are somewhat different for even and odd $n$. Taking into account that all the coefficients $\theta_n$ are zero for even $n$, we can rewrite recurrent relation 
(\ref{eq024}) as follows: 
\begin{equation}\label{eq027} 
\gamma_{2n} = \lambda \sum_{k=0}^{n-1} \gamma_{2n-2k-1} \; \theta_{2k+1} , \qquad n\ge 1, 
\end{equation}
\begin{equation}\label{eq0027} 
\gamma_1=1, \qquad \gamma_{2n+1} = \theta_{2n+1} +\lambda \sum_{k=0}^{n-1} \gamma_{2n-2k} \; \theta_{2k+1} , \qquad n\ge 1, 
\end{equation}
where, according to (\ref{eq0090}), 
\begin{equation}\label{eq00027}
\theta_1 = 1, \qquad  \theta_{2k+1} = (-1)^k \; \frac{\left( \frac{1}{2} \right)_k}{\left( \frac{m}{2} 
\right)_k} \; (c\Vert\boldsymbol\alpha\Vert)^{2k} , \qquad k\ge 1. 
\end{equation}

\section{Mixed moments} 

Series representations of the characteristic function $H(\boldsymbol\alpha,t)$ given above by Theorems 1 
and 2 enable us to calculate the mixed moments. Let $\bold q=(q_1,\dots, q_m)$ be a multi-index. Since at arbitrary fixed time $t>0$ the Markov random flight $\bold X(t) = (X_1(t),\dots,X_m(t))$ is concentrated 
in the ball $\bold B_{ct}^m$, this process, as well as its coordinates $X_j(t), \; j=1,\dots,m$, 
are bounded and, therefore, for arbitrary positive integer $k$, the condition $\Bbb E |X_j(t)|^k<\infty$ 
fulfills for all $j=1,\dots,m$. Then, for arbitrary fixed $t>0$, there exist the mixed moments   
$$\mu_{\bold q}(t) = \mu_{(q_1,\dots,q_m)}(t) = \Bbb E X_1^{q_1}(t)\dots X_m^{q_m}(t) , 
\qquad q_j\ge 0, \quad q_1+\dots +q_m \le k,$$ 
of the $m$-dimensional symmetric Markov random flight $\bold X(t)=(X_1(t),\dots,X_m(t))$ given by the formula: 
\begin{equation}\label{eq019} 
\mu_{\bold q}(t) = \mu_{(q_1,\dots,q_m)}(t) = (-i)^{q_1+\dots +q_m} \; 
\frac{\partial^{q_1+\dots +q_m}}{\partial\alpha_1^{q_1} \dots \partial\alpha_m^{q_m}} \; 
H(\alpha_1,\dots,\alpha_m; t) \biggr|_{\alpha_1= \dots =\alpha_m=0} . 
\end{equation}

Note, that if the characteristic function depends only on the square of the norm 
$\Vert\boldsymbol\alpha\Vert^2=\alpha_1^2+ \dots +\alpha_m^2$ of inversion parameter $\boldsymbol\alpha$, then, for arbitrary dimension $m\ge 2$ and arbitrary positive integer $n\ge 1$, the following useful relation holds: 
\begin{equation}\label{eq28} 
\frac{\partial^{2m}}{\partial\alpha_1^2 \dots \partial\alpha_m^2} \; 
(\alpha_1^2+ \dots +\alpha_m^2)^n \; \biggr|_{\alpha_1= \dots =\alpha_m=0} 
= \left\{ \aligned 
(2m)!! , \qquad & \text{if} \;\; n=m ,\\ 
0, \qquad & \text{if} \;\; n\neq m , \endaligned \right. \qquad m\ge 2, \quad n\ge 1. 
\end{equation}

Since the moment functions of the $2D$ and $4D$ Markov random flights are known (see \cite{kol03}), 
we will focus on evaluating the first and second mixed moments of the three-dimensional symmetric 
Markov random flight $\bold X(t)=(X_1(t),X_2(t),X_3(t)), \; t>0$. In this case $m=3$ and 
formula (\ref{eq28}) becomes: 
\begin{equation}\label{eq29} 
\frac{\partial^6}{\partial\alpha_1^2 \; \partial\alpha_2^2 \; \partial\alpha_3^2} \; 
(\alpha_1^2+\alpha_2^2+\alpha_3^2)^n \; \biggr|_{\alpha_1=\alpha_2=\alpha_3=0} 
= \left\{ \aligned 
6!! , \qquad & \text{if} \;\; n=3 ,\\ 
0, \qquad & \text{if} \;\; n\neq 3 , \endaligned \right. \qquad n\ge 1. 
\end{equation} 

To find these mixed moments, it is convenient to use the series representation (\ref{eq20}) of 
the characteristic function $H(t)$ rather than more complicated one (\ref{eq5}). Since each 
coefficient $\gamma_n=\gamma_n(\Vert\boldsymbol\alpha\Vert), \; n\ge 1,$ is some polynomial with respect 
to even powers of the variable $\Vert\boldsymbol\alpha\Vert$, then, taking into account that 
$$\frac{\partial^3}{\partial\alpha_1 \; \partial\alpha_2 \; \partial\alpha_3} \; 
\Vert\boldsymbol\alpha\Vert^{2k} \; \biggr|_{\alpha_1=\alpha_2=\alpha_3=0} = 
\frac{\partial^3}{\partial\alpha_1 \; \partial\alpha_2 \; \partial\alpha_3} \; 
(\alpha_1^2+\alpha_2^2+\alpha_3^2)^k \; \biggr|_{\alpha_1=\alpha_2=\alpha_3=0} = 0 , \qquad k\ge 1,$$
and   
$$\frac{\partial^3}{\partial\alpha_1 \; \partial\alpha_2 \; \partial\alpha_3} \; 
\gamma_n(\Vert\boldsymbol\alpha\Vert) \; \biggr|_{\alpha_1=\alpha_2=\alpha_3=0} =0 , \qquad n\ge 1,$$
we see that 
$$\frac{\partial^3}{\partial\alpha_1 \; \partial\alpha_2 \; \partial\alpha_3} \; 
H(\boldsymbol\alpha,t) \; \biggr|_{\alpha_1=\alpha_2=\alpha_3=0} =  
e^{-\lambda t} \sum_{n=0}^{\infty} 
\frac{\partial^3}{\partial\alpha_1 \; \partial\alpha_2 \; \partial\alpha_3} \; 
\gamma_{n+1}(\Vert\boldsymbol\alpha\Vert) \; \biggr|_{\alpha_1=\alpha_2=\alpha_3=0}  \; 
\frac{t^n}{n!} = 0 ,$$
and, therefore, for arbitrary fixed $t>0$, the first mixed moment corresponding to the multi-index 
$\bold q_1 = (1,1,1)$ is zero, that is, 
\begin{equation}\label{eq30} 
\mu_{\bold q_1}(t) = \mu_{(1,1,1)}(t) = 0 . 
\end{equation}

Let us now concentrate on evaluating the second mixed moment function 
$\mu_{\bold q_2}(t) = \mu_{(2,2,2)}(t)$ corresponding to the multi-index $\bold q_2 = (2,2,2)$. 
For the sake of brevity, let us denote the operator: 
\begin{equation}\label{eq31} 
\mathcal D = \frac{\partial^6 \; \cdot}{\partial\alpha_1^2 \; \partial\alpha_2^2 \; 
\partial\alpha_3^2} \; \biggr|_{\alpha_1=\alpha_2=\alpha_3=0}   
\end{equation}
In order to evaluate the second mixed moment, one needs to calculate the action of this operator on the coefficients $\gamma_n(\Vert\boldsymbol\alpha\Vert)$, that is, 
$$\mathcal D \; \gamma_n(\Vert\boldsymbol\alpha\Vert) = \frac{\partial^6 \; }{\partial\alpha_1^2 \; \partial\alpha_2^2 \; \partial\alpha_3^2} \; \gamma_n(\Vert\boldsymbol\alpha\Vert) \; 
\biggr|_{\alpha_1=\alpha_2=\alpha_3=0} , \qquad n\ge 2.$$
However, this is a very difficult analytical problem because we need an explicit form of the coefficients 
$\gamma_n(\Vert\boldsymbol\alpha\Vert), \; n\ge 2$, while they are given only in the determinant form 
(\ref{eq27}). Application of operator (\ref{eq31}) to such a determinant leads to an extremely complicated and cumbersome expression containing $n^6$ sub-determinants that, obviously, makes the problem impracticable (at least, analytically). 

That is why we confine ourselves to the case of several terms of series (\ref{eq20}). This will give us not exact, but asymptotic relation for the second mixed moment function $\mu_{(2,2,2)}(t)$ with respect to time variable. According to (\ref{eq20}), we have: 
\begin{equation}\label{eq32} 
\aligned 
H(\boldsymbol\alpha,t) = e^{-\lambda t} & \biggl[ \gamma_1(\boldsymbol\alpha) + 
\gamma_2(\boldsymbol\alpha) \; \frac{t}{1!}  + \gamma_3(\boldsymbol\alpha) \; \frac{t^2}{2!} + 
\gamma_4(\boldsymbol\alpha) \; \frac{t^3}{3!} + \gamma_5(\boldsymbol\alpha) \; \frac{t^4}{4!} \\ 
& \qquad\qquad\quad + \gamma_6(\boldsymbol\alpha) \; \frac{t^5}{5!} + \gamma_7(\boldsymbol\alpha) \; 
\frac{t^6}{6!} + \gamma_8(\boldsymbol\alpha) \; \frac{t^7}{7!}  + o(t^7) \biggr] .
\endaligned
\end{equation}

From (\ref{eq25}), for $m=3$, we get:
\begin{equation}\label{eq33} 
\aligned 
& \gamma_1(\boldsymbol\alpha) = 1, \qquad \gamma_2(\boldsymbol\alpha) = \lambda, \qquad 
\gamma_3(\boldsymbol\alpha) = \lambda^2 - \frac{1}{3} \; (c\Vert\boldsymbol\alpha\Vert)^2 , \\  
& \gamma_4(\boldsymbol\alpha) = \lambda^3 - \lambda \; \frac{2}{3} \; (c\Vert\boldsymbol\alpha\Vert)^2 , \qquad \gamma_5(\boldsymbol\alpha) = \lambda^4 - \lambda^2 \; (c\Vert\boldsymbol\alpha\Vert)^2 
+ \frac{1}{5} \; (c\Vert\boldsymbol\alpha\Vert)^4 \\
& \gamma_6(\boldsymbol\alpha) = \lambda^5 - \lambda^3 \; \frac{4}{3} \; (c\Vert\boldsymbol\alpha\Vert)^2 
+ \lambda \; \frac{23}{45} \; (c\Vert\boldsymbol\alpha\Vert)^4 \\ 
& \gamma_7(\boldsymbol\alpha) = \lambda^6 - \lambda^4 \; \frac{5}{3} \; (c\Vert\boldsymbol\alpha\Vert)^2 
+ \lambda^2 \; \frac{14}{15} \; (c\Vert\boldsymbol\alpha\Vert)^4  
- \frac{1}{7} \; (c\Vert\boldsymbol\alpha\Vert)^6 \\ 
& \gamma_8(\boldsymbol\alpha) = \lambda^7 - 2\lambda^5 \; (c\Vert\boldsymbol\alpha\Vert)^2  
+ \lambda^3 \; \frac{22}{15} \; (c\Vert\boldsymbol\alpha\Vert)^4  
- \lambda \; \frac{44}{105} \; (c\Vert\boldsymbol\alpha\Vert)^6  
\endaligned 
\end{equation}
In view of (\ref{eq29}), we have: 
$$\mathcal D \; \gamma_1(\boldsymbol\alpha) = \mathcal D \; \gamma_2(\boldsymbol\alpha) = \mathcal D \; \gamma_3(\boldsymbol\alpha) = \mathcal D \; \gamma_4(\boldsymbol\alpha) = \mathcal D \; 
\gamma_5(\boldsymbol\alpha) = \mathcal D \; \gamma_6(\boldsymbol\alpha) = 0 ,$$
while 
$$\mathcal D \; \gamma_7(\boldsymbol\alpha) = - \frac{48}{7} \; c^6 , \qquad 
\mathcal D \; \gamma_8(\boldsymbol\alpha) = - \lambda \; \frac{704}{35} \; c^6 .$$

Therefore, applying operator $\mathcal D$ to (\ref{eq32}), we obtain: 
$$\aligned 
\mathcal D \; H(\boldsymbol\alpha,t) & = - e^{-\lambda t} \biggl[ \frac{48}{7} \; c^6 \; \frac{t^6}{6!} 
+ \lambda \; \frac{704}{35} \; c^6 \; \frac{t^7}{7!} + o(t^7) \biggr] \\
& = - e^{-\lambda t} \biggl[ \frac{1}{105} \; (ct)^6 + \lambda \; \frac{44}{11025} \; c^6 t^7 
+ o(t^7) \biggr] .
\endaligned$$
Then, according to (\ref{eq019}), and taking into account that $i^6 = -1$, we finally get: 
\begin{equation}\label{eq34}
\mu_{\bold q_2}(t) = \mu_{(2,2,2)}(t) = e^{-\lambda t} \biggl[ \frac{1}{105} \; (ct)^6 + \lambda \; 
\frac{44}{11025} \; c^6 t^7 + o(t^7) \biggr] .
\end{equation}

Thus, we can summarize the above results (\ref{eq30}) and (\ref{eq34}) in the following theorem. 

\bigskip

{\bf Theorem 3.} {\it For arbitrary $t>0$, the first mixed moment of the three-dimensional symmetric 
Markov random flight $\bold X(t) = (X_1(t),X_2(t),X_3(t))$ is identically equal to zero, that is,} 
\begin{equation}\label{eq35} 
\mu_{\bold q_1}(t) = \mu_{(1,1,1)}(t) = 0 . 
\end{equation}

{\it For the second mixed moment function of $\bold X(t)$, the following asymptotic relation holds:}
\begin{equation}\label{eq36}
\mu_{\bold q_2}(t) = \mu_{(2,2,2)}(t) = e^{-\lambda t} \biggl[ \frac{1}{105} \; (ct)^6 + \lambda \; 
\frac{44}{11025} \; c^6 t^7 + o(t^7) \biggr] .
\end{equation} 

\bigskip 

{\bf Remark 7.} Asymptotic formula (\ref{eq36}) provides a very good approximation of the second mixed 
moment, especially for small values of time $t$. For example, for the time value $t=0.3$, relation 
(\ref{eq36}) yields (for $c=3, \; \lambda=5$) the estimate: $\mu_{(2,2,2)}(0.3) \approx 0.00183921$ and 
the error of this estimate does not exceed the value 0.00001464 multiplied by some constant.


\begin{thebibliography}{18} 

\bibitem{bras}
Brasiello A., Crescitelli S., Giona M. One-dimensional hyperbolic transport: 
positivity and admissible boundary conditions derived from the wave formulation. 
{\it Phys. A}, {\bf 449} (2016), 176-191.

\bibitem{broad1}
Broadbridge P., Kolesnik A.D., Leonenko N., Olenko A., Omari D. Spherically restricted 
random hyperbolic diffusion. {\it Entropy}, {\bf 22} (2020), 217-248.   

\bibitem{broad2}
Broadbridge P., Kolesnik A.D., Leonenko N., Olenko A. Random spherical 
hyperbolic diffusion. {\it J. Statist. Phys.}, {\bf 177} (2019), 889-916.  

\bibitem{cane1}
Cane V. Random walks and physical processes. {\it Bull. Int. Statist.
Inst.}, {\bf 42} (1967), 622-640.

\bibitem{cane2}
Cane V. Diffusion models with relativity effects. In: {\it
Perspectives in Probability and Statistics}, Applied Probability
Trust, 1975, Sheffield, pp. 263-273. 

\bibitem{eon} 
d'Eon E., McCormick N.J. Radioactive transfer in half spaces of arbitrary dimensions. 
{\it J. Comput. Theoret. Transport}, {\bf 48} (2019), 280-337. 

\bibitem{giona1}
Giona M., Brasiello A., Crescitelli S. Markovian nature, completeness, regularity and
correlation properties of generalized Poisson-Kac processes. {\it J. Stat. Mech.: Theory and Experiment},
{\bf 2} (2017), iss. 2, 023205.

\bibitem{giona2}
Giona M., Brasiello A., Crescitelli S. Stochastic foundations of undulatory transport phenomena:
generalized Poisson-Kac processes - Part I: Basic theory. {\it J. Phys. A: Mathematical and Theoretical},
{\bf 50} (2017), iss. 33, 335002.

\bibitem{giona3}
Giona M., Brasiello A., Crescitelli S. Stochastic foundations of undulatory transport phenomena:
generalized Poisson-Kac processes - Part II: Irreversibility, norms and entropies.
{\it J. Phys. A: Mathematical and Theoretical}, {\bf 50} (2017), iss. 33, 335003.

\bibitem{giona4}
Giona M., Brasiello A., Crescitelli S. Stochastic foundations of undulatory transport phenomena:
generalized Poisson-Kac processes - Part III: Extensions and applications to kinetic theory and transport.
{\it J. Phys. A: Mathematical and Theoretical}, {\bf 50} (2017), iss. 33, 335004.

\bibitem{gold}
Goldstein S. On diffusion by discontinuous movements and on the
telegraph equation. {\it Quart. J. Mech. Appl. Math.}, {\bf 4} (1951), 129-156. 

\bibitem{gr}
Gradshteyn I.S., Ryzhik I.M. {\it Tables of
Integrals, Series and Products.} Academic Press, 1980, NY. 

\bibitem{kac}
Kac M. A stochastic model related to the telegrapher's equation.
{\it Rocky Mount. J. Math.}, {\bf 4} (1974), 497-509. (Reprinted from: 
Kac M. Some stochastic problems in physics and mathematics. 
In: {\it Magnolia Petroleum Company Colloquium Lectures
in the Pure and Applied Sciences}, No. 2, October 1956). 

\bibitem{kol1}
Kolesnik A.D. {\it Markov Random Flights.} Taylor \& Francis Group/CRC Press, 2021, 
London-New York-Boca Raton.

\bibitem{kol2}
Kolesnik A.D. The explicit probability distribution of a six-dimensional random flight. 
{\it Theory Stoch. Process.,} {\bf 15(31)} (2009), 33-39.

\bibitem{kol3}
Kolesnik A.D. Random motions at finite speed in higher dimensions. 
{\it J. Statist. Phys.,} {\bf 131} (2008), 1039-1065. 

\bibitem{kol03}
Kolesnik A.D. Moments of the Markovian random evolutions in two and four dimensions. 
{\it Bull. Acad. Sci. Moldova, Ser. Math.}, {\bf 2(57)} (2008), 68-80.  

\bibitem{kol4}
Kolesnik A.D. A four-dimensional random motion at finite speed. 
{\it J. Appl. Probab.,} {\bf 43} (2006), 1107-1118. 

\bibitem{kol5}
Kolesnik A.D., Orsingher E. A planar random motion with 
an infinite number of directions controlled by the damped wave equation. 
{\it J. Appl. Probab.}, {\bf 42} (2005), 1168-1182. 

\bibitem{kolrat}
Kolesnik A.D., Ratanov N. {\it Telegraph Processes and Option Pricing.} Springer, 2013, 
Berlin-Heidelberg (see also the 2nd edition, Springer, 2022). 

\bibitem{korn}
Korn G.A., Korn T.M. {\it Mathematical Handbook}. McGraw-Hill, 1968, NY. 

\bibitem{mas}
Masoliver J., Porr\'a J.M., Weiss G.H. 
Some two and three-dimensional persistent random walks. 
{\it Physica A}, {\bf 193} (1993), 469-482. 

\bibitem{pbm1} 
Prudnikov A.P., Brychkov Yu.A., Marichev O.I. {\it Integrals and Series. Supplementary Chapters.} 
Nauka, 1986, Moscow. (In Russian)

\bibitem{sta1}
Stadje W. Exact probability distributions for non-correlated  
random walk models. {\it J. Statist. Phys.}, {\bf 56} (1989), 415-435.

\bibitem{sta2}
Stadje W. The exact probability distribution of a two-dimensional 
random walk. {\it J. Statist. Phys.}, {\bf 46} (1987), 207-216. 


\end{thebibliography}
\end{document}